\documentclass[russian,english]{article}
\usepackage[T1]{fontenc}
\usepackage[koi8-r]{inputenc}
\usepackage{amsmath}
\usepackage{graphicx}
\usepackage{amssymb}
\usepackage{esint}

\makeatletter
\newtheorem{theorem}{Theorem}

\newtheorem{exer}{Exersize}
\newtheorem{lemma}{Lemma}

\newtheorem{problem}{Problem}

\newtheorem{remark}{Remark}

\makeatother

\makeatother

\makeatother

\usepackage{babel}

\makeatother

\makeatother

\usepackage{babel}

\makeatother

\usepackage{babel}

\makeatother
\usepackage{babel}\date{}

\makeatother

\usepackage{babel}

\makeatother

\usepackage{babel}

\begin{document}

\title{Introduction to stochastic models of transportation flows. Part I. }

\author{V. A. Malyshev, A. A. Zamyatin}
\maketitle
\begin{abstract}
We consider here probabilistic models of transportation flows. The
main goal of this introduction is rather not to present various techniques
for problem solving but  to present some intuition to invent adequate
and natural models having visual simplicity and simple (but rigorous)
formulation, the main objects being cars not abstract flows. The papers
consists of three parts. First part considers mainly linear flows
on short time scale - dynamics of the flow due to changing driver
behavior. Second part studies linear flow on longer time scales -
individual car trajectory from entry to exit from the road. Part three
considers collective car movement in complex transport networks.
\end{abstract}
\tableofcontents{}

\bigskip{}

Mathematical models of car traffic can be very different: from partial
differential equations to modern computer graphics where points move
on video along the edges of some graph. We consider here only probabilistic
models. The main goal of this introduction is rather not to present
various techniques for problem solving but  to present some intuition
to invent adequate and natural models having visual simplicity and
simple (but rigorous) formulation, the main objects being cars not
abstract flows. Moreover, the postulates and parameters in these models
should allow statistical verification and estimation, and should not
use doubtful physical analogies. Probabilistic models should be tightly
related to psychology of drivers, if the drivers are not robots. There
is no such theory, and we present an attempt to start it.

Probabilistic approach to transportation problems exists already more
than 50 years, see \cite{Haight,Renyi,SolomonWang}, however in this
paper we follow newer probability theories and consider more difficult
problems. There are however some works which, for various reasons,
we did not include in our introduction, see for example \cite{Blank,MaxPlus}

The papers consists of three parts. First part considers mainly linear
flows on short time scale - dynamics of the flow due to changing driver
behavior. Second part studies linear flow on longer time scales -
individual car trajectory from entry to exit from the road. Part three
considers collective car movement in complex transport networks.

\section{Car flows}

\subsection{Marked point fields}

Current is normally the mean number $J$ of cars crossing, per unit
of time and in the given direction, some section of the road. Flow
is normally either a static random configuration\[
...<x_{i}<x_{i-1}<...\]
 of cars at a given time moment, or a probability measure on the set
of car trajectories $\{x_{i}(t)\}$.

\paragraph{What is a car configuration}

In this paper the maximally detailed description for the configuration
of cars at a given time moment will be as follows. The car is specified
by some index $\alpha$. For example, a road (route) with $k$ (traffic)
lanes $1,2,...,k$ , is represented by $k$ lines, parallel to $x$-axis.
Then the index $\alpha=(m,i)$ marks out the $i$-th car on the lane
$m$. Index $i$ enumerates cars on the lane so that $i$-th car follows
car with index $i-1$. Let $d_{\alpha}$ the length of the car indexed
by $\alpha$, $x_{\alpha}(t)$ - its coordinate (for example, of the
front buffer). The cars move in the positive $x$-direction. Further
on we omit the lane index (the reader can add it whenever necessary)
and use only index $i$.

Denote \[
d_{i}^{+}(t)=x_{i-1}(t)-x_{i}(t)-d_{i-1}\]
the distance of the car $i$ to the previous car at time $t$. The
distance of the car to the following car \[
d_{i}^{-}(t)=d_{i+1}^{+}(t)\]
is also important for the driver.

\paragraph{How probabilities on the car configurations are introduced}

Formally, the point flow on the real line $R$ is given by the probability
measure on the set of all countable locally finite (that is finite
on any bounded interval) subsets of $R$. Otherwise speaking, it is
given by compatible systems of probabilities\[
P(I_{1},k_{1};...;I_{n},k_{n})\]
 of the events that in the intervals $I_{j},j=1,...,n$ there are
exactly $k_{j}$ particles.

The main question of course is how to specify this system more concretely.
There are two big parts of probability theory which use different
ways of doing this. These are renewal theory (see for example, \cite{KoksSmit})
and Gibbs point fields theory \cite{Ruelle,MalMin}. The first is
essentially simpler but works only in one-dimensional case. The second
has deep relations with physics, is suitable also for multidimensional
point configurations but is more complicated, and we shall not touch
it here.

The simplest point flow is the Poisson flow, see \cite{Kingman}.
Most natural way to understand it is as follows. Consider the interval
$[-N,N]$ and throw (choose) on it independently and randomly (more
exactly, uniformly) $M=[\rho N]$ points, where $\rho>0$ is some
constant called density. It is easy to calculate the binomial probability
$P_{N,M}(k,I)$ that exactly $k$ points will be inside the finite
interval $I$. This probability tends, as $N\to\infty$, to the Poisson
expression (Poisson distribution for given $I$)\[
P(k,I)=\frac{\{\rho|I|\}^{k}}{k!}e^{-\rho|I|}.\]
 More general flows can be easily constructed on the half line $[0,\infty)$.
Namely, random points\[
x_{0}=0,x_{1},...,x_{n},...\]
 are defined as the sums \[
x_{1}=\xi_{1},x_{2}=\xi_{1}+\xi_{2},....\]
 of i.i.d. (independent identically distributed) random variables
$\xi_{i}>0,i=1,2,...$, with the distribution function $G(x)=Pr(\xi_{i}<x)$.

To define a translation invariant point flow on all real line one
problem remains - how to choose point $0$, from which one can consecutively
place independent variables to the right and to the left. For this
one should use the following (in fact, one of the main) assertion
of the renewal theory. Let $P(t,t+\Delta t)$ be the probability that
the interval $(t,t+\Delta t)$ will contain exactly one point. Then,
if $\xi_{i}$ have density, then on the half-line the limit $\lim_{\Delta t\to0}\frac{1}{\Delta t}P(t,t+\Delta t)$
exists and tends (as $t\to\infty$) to $\rho=(E\xi)^{-1}$. Then the
limit (as $t\to\infty$) probability density of the event $A(t,s)$,
that the distance from point $t$ to the first random point $x_{i}>t$
is greater than some $s>0$, equals\begin{equation}
\rho Pr(\xi_{i}=x_{i}-x_{i+1}>s)=\rho(1-G(s)).\label{renewaldensity}\end{equation}
 That is why the first (after $0$) point $x_{1}$ should be taken
on the random distance from $0$ with density (\ref{renewaldensity}).
The following consecutive distances should be independent with distribution
$G(x)$.

\paragraph{Alternating flows}

The distances in-between neighbor points are not necessarily identically
distributed. The distributions can alternate. Take, for example, two
sequences of i.i.d. random variables $\xi_{1},\xi_{2},...$ and $\eta_{1},\eta_{2},...$,
and put\[
x_{2n}=\xi_{1}+...+\xi_{n}+\eta_{1}+...+\eta_{n},x_{2n-1}=\xi_{1}+...+\xi_{n}+\eta_{1}+...+\eta_{n-1}.\]
 Then the construction of the flow on the whole line can be done as
above.

In our case the (i.i.d.) lengths of cars $d_{i}$ and (i.i.d.) functions
$d_{i}^{+}$ alternate.

\paragraph{Marked flows}

To each point $x_{i}$ of the point process one can assign some variable
$\sigma_{i}$, taking values in some set $S$. This variable is called
a mark or a spin of the point (particle, car) $i$, then one can call
the new object random marked point set (flow, process). It is defined
by some probability measure on the set of sequences of pairs $(x_{i},\sigma_{i})$.
Often the measure on the countable sets is given, that is the flow
without marks being given, then $\sigma_{i}$ are considered to be
i.i.d. random variables. In section 1.5 we shall consider marked process
where the marks will be the velocities of the cars, moreover with
sufficiently complicated correlations of the velocities with the coordinates.

\paragraph{About Markov processes}

The models of time evolution of cars often use Markov processes and
it is useful to remind some terminology. However, the mere definitions
of the Markov process and its properties (for example, ergodicity)
may differ. We give the definition we use for discrete time.

Consider on some phase space $X$ a system of measures (transition
probabilities) $P(A|x)$, meaning the probability that the process
at time $t+1$ hits the set $A\subset X$, if at time $t$ the process
was at the state $x\in X$. If all measures $P(A|x)$ are one-point,
then it is equivalent to the deterministic map $T:X\to X$, that is
$P(.|x)$ is the unit measure at the point $T(x)$. This deterministic
map defines a dynamical system.

Note that the system of measures $P(A|x)$ defines the mapping $U$
of the set of all measures on $X$ to itself\[
U\mu=\int P(.|x)d\mu(x)\]
 The notion of invariant (with respect to $U$) measure is very important.
Normally one studies its existence, uniqueness and other properties.

Using the system of transition probabilities one can construct various
sequences of random variables $\xi_{n}=x_{n}$ with values in $X$,
or their distributions $\mu_{n}$ on $X$, where $\Pr(\xi_{n}\in A)=\mu_{n}(A)$.
The probability space is the set of trajectories $\{x_{n}\}$. For
example, it can be stationary Markov process $\xi_{n},n\in Z,$ where
$\mu_{n}$ are invariant measures on $X$, or non-stationary process
$\xi_{n},n\in Z_{+}$, with given initial distribution $\mu_{0}$
on $X$.

Markov process can be understood either as ONE of this sequences $\xi_{n}$,
or as the whole family of such sequences. Correspondingly, the terminology
(for example the notion of ergodicity) differs. Dynamical system with
given invariant measure $\mu$ is called ergodic if any invariant
subset has $\mu$ zero or one. Any stationary process can be considered
as the dynamical system - time shift on the set of trajectories. Then
the ergodicity can be defined as the ergodicity of this dynamical
system. However more often, Markov process is understood more generally
and the notion of ergodicity is defined differently. Namely, the process
is called ergodic if there exists the only invariant measure $\mu$
on $X$, and for any initial measure $\mu_{0}$ as $n\to\infty$ we
have $\mu_{n}\to\mu$ in the sense of weak convergence.

Note that Markov processes normally belong to one of two classes.
The first class includes all processes such that there exists a positive
measure on $X$ (may be infinite), with respect to which all measures
$P(.|x)$ are absolutely continuous. This includes all classical Markov
process - finite and countable Markov chains, diffusion processes,
etc. Such process are ergodic if there are no non-trivial invariant
subsets (this property is called irreducibility), and there is a unique
invariant measure. In most cases convergence to it from any initial
state follows. For countable chains an equivalent condition is positive
recurrence that is finiteness (for any pair $x,y\in X$) of the mean
hitting time of $x$ starting from $y$.

The second class is characterized by the property that all measures
$P(.|x)$ are mutually singular. All infinite particle processes are
in this class. The theory of such processes is essentially more complicated.

\subsection{Traffic capacity as a function of velocity }

\paragraph{Driver's psychology in the simplest flow}

Detailed modeling of the driver psychology is impossible but some
relation is obvious. Thus, the driver $i$ sees several cars (often
only one car in front of him) in the flow and chooses subjectively
optimal distance to the preceding car. If the velocity $v_{i-1}(t)=\frac{dx_{i-1}(t)}{dt}$
of the preceding car changes slowly, then one can assume that the
driver's reaction is faster, and the chosen distance $D_{i}^{+}$
depends of the velocity of the preceding car at the same time moment\[
D_{i}^{+}=D_{i}^{+}(v_{i-1})\]
 (index $i$ says that the functions $D_{i}^{+}(v)$ can be different
for different drivers). Let us call the flow algorithmic at time $t$,
if for all $i$\[
d_{i}^{+}(t)=D_{i}^{+}(v_{i-1}(t)),\]
 that is all velocities are consecutively defined by the velocity
of the previous cars. Of course, individual functions can be interesting
for policemen, but we are interested with their statistical characteristics.
In probabilistic approach the functions $D_{i}^{+}(v)$ become i.i.d.
random variables, depending on the velocity $v$ of the previous car
as a parameter. The distribution of these functions cannot be deduced
from mathematical or physical laws, but rather it depends on the individual
and collective psychology of drivers and should be found experimentally,
see \cite{Revised}.

\paragraph{Deterministic dynamics without overtaking}

If all cars, drivers and velocities $v$ are identical, then many
problems appear to be very simple. Let $d$ be the length of cars
and $d^{+}=D^{+}(v)$ be the distance to the previous car, which the
driver strictly follows. Already such dynamics allows to understand
many qualitative aspects.

Let us define the road capacity as maximally possible current through
it\[
J_{max}=\max_{v}v\lambda(v),\]
 where maximum is over the allowed interval of velocities, and where\[
\lambda(v)=\frac{k}{d+D^{+}(v)}\]
 is the density of cars on $k$-lane road for a fixed velocity $v$.
It is clear that the flow capacity can become smaller when the velocity
increases. This simple conclusion says only that many drivers increase
the distance to the previous car when the velocity becomes too big.

\paragraph{Random dynamics without overtaking}

We get the same conclusion if the velocities $v$ are the same, but
the functions $d_{i}^{+}$ are random independent, and their expectations
(for given $v$) are equal to some number $d^{+}(v)$. We see that
the fact of non-trivial dependence of the road capacity from the velocity
is in fact the trivial consequence of the driver behavior and does
not need any models. However, for nicer questions stochastic models
are necessary. Now we will introduce a sufficiently general probabilistic
model with a rich spectrum of phases. The well-known exclusion processes
appear as a degenerate particular case. Other models see in \cite{Blank,Tanaka,Revised}.

\subsection{Random dynamics with overtaking (random grammars)}

Here we have an interesting connection with recently discovered object
- random grammars, see \cite{Grammars}. We give short substantial
description of one such model.

Let at time $t=0$ all cars are situated on the left half-axis, and
the movement is one-lane. We subdivide the lane onto intervals (cells,
enumerated by $Z$) of fixed length and assume that there can be at
most one car in the cell. Thus, finite sequence (group, cluster) of
cars is identified with the pair $(S,r)$, where $r\in Z$, and $S$
is a finite sequence (word) of three symbols $0,1,2$\[
S=s_{N}...s_{2}s_{1}.\]
 Here $0$ corresponds to empty cell, $1$ - to active (fast) driver
in the cell, $2$ - to a quiet driver. The length $N=N(t)$ of the
word and all symbols $s_{k}(t)$  change with time, with the restriction
that always $s_{1}(t)\neq0$ for all $t\geq0$. At any time $t$ any
symbol $s_{k}(t)$ has coordinate $x(s_{k}(t))$. The coordinates
\begin{equation}
x(s_{k}(t))=x(s_{1}(t))-k+1\label{coordinates}\end{equation}
 are uniquely defined by the coordinate $x(s_{1}(t))$ of the first
symbol, denoted by $r=r(t)$.

This dynamics models the process of acceleration and slowdown of different
drivers. It is defined as continuous time Markov chain $(S(t),r(t))$
on the set $\{(S,r)\}$ of pairs. Jump intensities are defined as
follows. Process $r(t)$ models the movement of all group with common
velocity $v$. Namely, $r$ is increased by one with probability $vdt$
for time $dt$, and all coordinates immediately change correspondingly
to formula (\ref{coordinates}). Thus, $S(t)$ describes replacements
relative to some common motion, and is given by some random grammar,
that is by the list of possible local substitutions (5 substitution
types), that is a subword of $S$ to another subword. Any substitutions
of this list are independent and have different intensities (4 parameters).
Here is the list: 
\begin{enumerate}
\item $10\to01$ - fast driver moves to empty place in front of him, simultaneously
creating additional free place after him, with probability $\lambda_{0}^{+}dt$
for time $dt$; 
\item $120\to021$ - fast driver overtakes quiet driver with intensity $\lambda_{1}^{+}dt$; 
\item $22\to202$, $21\to201$ - a cautious driver brakes, increasing the
distance in front of him with probability $\lambda_{2}^{-}dt$. Note
that the length $S$ increases here (extra empty cell appears), that
causes the one-cell shift of all cars behind. This is a non-local
jump, in fact such rearrangement takes some time, but we assumed that
this time is on a shorter scale; 
\item $200\to020$ - a quiet driver accelerates with probability $\lambda_{2}^{+}dt$
(if he decides that in front of him there is too much free place). 
\end{enumerate}
It is necessary to say that for exact formulation of results which
we only shortly describe and to get various qualitative types (phases)
of movement, one needs to do various scalings of parameters $t,N,\lambda$,
We shortly describe only 3 phases. 
\begin{enumerate}
\item If $\lambda_{2}^{\pm}$ are small relative to other two parameters,
then cars of type 2 move synchronously with common velocity, and fast
cars have additional relative velocity. If the number of fast drivers
is small, then this relative velocity is defined by the movement of
one car among static obstacles and depends on the density $\rho_{2}$
of type 2 cars and density of holes $\rho_{0}$, it is approximately
equal to\[
v_{rel}=\lambda_{0}^{+}\rho_{0}+2\lambda_{1}^{+}\rho_{2}\]

\item If $\lambda_{2}^{-}$ is small with respect to other parameters (no
non-local effects), and $\lambda_{2}^{+}$ has the same order as $\lambda_{0}^{+},\lambda_{1}^{+}$,
then the difference between types disappears. We have a process close
to the so called TASEP - totally asymmetric exclusion process. If
\[
\lambda_{0}^{+}=\lambda_{1}^{+},\lambda_{2}^{-}=0\]
 then it coincides with TASEP (about TASEP see \cite{Blank}). 
\item If $\lambda_{2}^{+}$ is small, and $\lambda_{2}^{-}$ is large with
respect to other two parameters, then the picture is different. Any
overtaking $120\to021$ immediately forces braking of the car $2$
and, as a consequence, ALL following cars slow down. For the cars
closer to the end of the word the slow down will be very essential,
if the flow is sufficiently dense (few cells with zeros), as many
type 2 cars will brake. 
\end{enumerate}
One can complicate the introduced dynamics, for example, to avoid
discretization (see the end of this section), introducing positive
real numbers instead of zeroes - distances between consecutive cars.
This demands essential reformulation (see  section 1.6), in particular
for type 3 jumps, but rough qualitative effects will be the same.

\subsection{Growth of a jam}

If the incoming traffic flow to some fixed domain equals $J_{in}$,
and outgoing flow is $J_{out}<J_{in}$, the number of cars in this
domain increases in time $t$ on \[
t(J_{in}-J_{out})\]
 This is true, however, only if this domain is not situated on the
road itself. For example, if the number of cars in the jam on the
road increases, the answer is different..The problem is that the domain
itself can grow because of the arriving cars. To make this more precise
one should find an appropriate model.

Assume that the cars of the same length $d$ move on one lane road
with velocity $v$ on the same distance $d^{+}$ between consecutive
cars. During some time $t$ the movement has been stopped by some
obstacle, for example by red traffic light. Assume that any car stops
on the distance $d_{0}^{+}<d^{+}$ from the previous car.

\begin{exer} Prove that during time $t\to\infty$ the jam (that is
maximal length $L(t)$ of the part of road, where all cars stand still)
before the obstacle will have the length asymptotically equal to\begin{equation}
L(t)\sim_{t\to\infty}tv\frac{d+d_{0}^{+}}{d^{+}-d_{0}^{+}}.\label{rostProbki}\end{equation}
 \end{exer}

It seems that this result depends only on the mean values and remains
true even if overtaking is possible. This was proved in \cite{MalZam}
for independent movement of the cars (that is when the cars do not
hinder each other), moreover the car velocities can fluctuate but
have the same values and equal $v$. The proof is absolutely not trivial.

\subparagraph{Local widenings and contractions of the road}

What occurs when a part of road with $k$ lanes turns to the part
with $l$ lanes. Let this occur at the point $x=0$.

Case $k<l$. Assume that maximally allowed velocity equals $v_{max}$
and the drivers are assumed to be disciplined. All cars move along
$k$-lane road with velocity $v<v_{max}$. And it is impossible to
move faster due to the fundamental relation \[
d+D^{+}(v)=\rho^{-1}\]
 between density $\rho$ of cars and their velocity. Then along $l$-lane
road the cars could conserve the density and move with the same speed,
but $\rho$ can change so that the cars will move with some greater
speed $v_{1}$. The time gain is\[
\frac{L}{v}-\frac{L}{v_{1}}\]
 Case $k>l$. Here three different situations are possible.

\textbf{Free flow} If the flow is very sparse, then the cars arrive
at the point $0$ alone and will not notice the decrease in the number
of lanes.

\textbf{Growing jam} Denote $J_{k}$ the incoming current and let
$J_{l,max}$ be maximally possible current along $l$-lane road. If
$J_{k}>J_{l,max}$, then the jam will grow, and the number of cars
in this jam will grow as $t(J_{k}-J_{l,max})$, more exactly as in
formula (\ref{rostProbki}).

\textbf{Delay} For the case $J_{k}<J_{l,max}$ practical observations
say that at point $0$ jams of random lengths can appear, which however
will not grow too much. There are no corresponding stochastic models.
To create such model, it is useful to have a collection of models
for shorter time phenomenon, Namely, when the standing cars at the
traffic lights start their movement. We shall describe some models
of this type now.

\subsection{Speeding up }

In \cite{FerrPech} the cars are points \[
...<x_{i}(t)<x_{i-1}(t)<...\]
 on the line. At initial time $t=0$ the cars stand still and are
described by Poisson point field with density $\rho<1$. The cars
can have two velocities: 0 and 1, and overtaking is prohibited. Any
incumbent car starts its movement with speed $1$, independently of
the others, in exponential time with mean 1. It can occur that the
car with number $i$ will reach the car $i-1$ when it did not start
its movement. Then the car $i$ stops and resumes its movement in
exponential time after the moment when the car $i-1$ starts movement.
This rule acts at any time. The model describes roughly how clusters
of cars start to move.

The main result is that with probability $1$ each car will stop finite
number of times (if $\rho<1$). Let $t_{i}$ be the minimal time moment
after which the car $i$ does not stop anymore. Then for any $i$
and $k$ and any time moments $t_{i},t_{i-1},...,t_{i-k}$ the random
variables\[
x_{i-1}(t_{i-1})-x_{i}(t_{i}),x_{i-2}(t_{i-2})-x_{i-1}(t_{i-1}),...,x_{i-k}(t_{i-k})-x_{i-k+1}(t_{i-k+1})\]
 will be independent and exponentially distributed. Otherwise speaking
after leaving the jam, the cars will form Poisson configuration with
the same density $\rho$.

Assume now that at time $0$ all cars are situated on the left half-axis.
The distribution is again Poisson with density $\rho$. Each point
moves with velocity $v>0$, if the distance to the next car from the
right is not less than some $d_{eff}>0$, and stands still otherwise.
Here it is evident that any particle will not stop starting from some
time moment. But in this model one can get more. Consider the following
random variables: random time $\tau_{k}^{(1)}$ when $k$-th point
starts movement, random time $\tau_{k}^{(2)}$ when this point does
not stop anymore, the distance $x_{k}$ between points $k$ and $1$,
starting from time moment $\tau_{k}^{(2)}$.

\begin{problem}

Find asymptotics of distributions of these random variables when $k\to\infty$.

\end{problem}

The connection with delay problem is evident. Let we have two lanes
and on any lane the flow density be $\rho$; the joint flow thus have
density $2\rho$. The cars from the first lane now want to squeeze
into the second. The squeezing algorithms can be quite different.
For example, any car can squeeze independently of the others, if its
distances (along $x$-axis) to the cars from the second lane are not
less than some $d^{+}$.

\subsection{Short and long range order with time dependent velocities}

Here the cars are represented by points $x_{i}$. To the car $i$
a stationary random process $w_{i}(t)$ is assigned, which defines
time evolution of its velocity on empty road (that is when there are
no obstacles in front of this car). This process implicitly defines
the activity of the driver at given time $t$. The processes $w_{i}(t)$
are mutually independent and are defined only by the psychology of
the individual driver. Assume that there are constants $0<C_{1}<C_{2}<\infty$
such that for all $t,i$ \[
C_{1}<w_{i}(t)<C_{2}.\]

The flow is given by initial coordinates $x_{i}(0)$ of the cars,
and their movement is defined as\[
x_{i}(t)=x_{i}(0)+\int_{0}^{t}v_{i}(s)ds,\]
 where $v_{i}(t)$ is the velocity of $i$-th car, defined below.
Moreover, we assume that the initial coordinates are such that the
distances $d_{i}^{+}(0)$ are independent and, for example, are exponentially
distributed with parameter $\rho(0)$.

We shall say that the car $i$ has an obstacle in front of itself
at time $t$, if \[
x_{i}(t-0)=x_{i-1}(t)\]
 The process will be completely defined if for any $t_{1},...,t_{n},i_{1},...,i_{n}$
we could define finite-dimensional distribution of the vectors\[
(v_{i_{1}}(t_{1}),...,v_{i_{n}}(t_{n}))\]
 where some $i_{k}$ can coincide. The following rules do this job: 
\begin{enumerate}
\item (rule for free road) if none of the cars $i_{1},...,i_{k}$ for $k\leq n$
does not have in front of itself an obstacle, then the distribution
of the vector $v_{i_{1}}(t_{1}),...,v_{i_{k}}(t_{k})$ coincides with
the distribution of the vector $w_{i_{1}}(t_{1}),...,w_{i_{k}}(t_{k})$
and is independent of the distribution of the vector $v_{i_{k+1}}(t_{k+1}),...,v_{i_{n}}(t_{n})$; 
\item (rule for an obstacle) if the car $i$ has obstacle in front of it
at time $t$, then $v_{i}(t)=v_{i-1}(t)$; 
\item (rule for overtaking) if the car $i$ has obstacle in front of it
at time $t$, then it overtakes this obstacle with some intensity
$\lambda$, during the whole (random) time interval while $w_{i}(t)>v_{i-1}(t)$.
The meaning of this condition is that the driver overtakes if his
activity is higher. 
\end{enumerate}
Even for such simplest behavior of drivers and corresponding definition
of the flow there are many problems. We formulate some of them.

Let us call \textbf{free phase} the case when the intensity of overtaking
equals infinity. Then for any cars $i,j$, their velocities are independent
and the covariance is zero, that is\[
Cov_{ij}(t)=Ev_{i}(t)v_{j}(t)-Ev_{i}(t)Ev_{j}(t)=0\]
 The \textbf{complete order} \textbf{phase} is when $\lambda=0$.
Then there is no overtaking at all and the cars move as a queue.

\begin{problem}

Assume that all $w_{i}(t)$ have the same distribution. Is it true
that for small $\lambda$ the movement will resemble complete order
phase, and the large $\lambda$ will resemble free phase ? Does exist
some third phase (cluster phase) for intermediate values of $\lambda$
?

\end{problem}

Possibly to get more interesting qualitative picture one has to define
the same process with the lengths $d_{i},d_{i}^{+}$, and with additional
indices corresponding to several lanes. Also the behavior of the driver
can depend not only on the car in front but also of the car just behind.

\subsection{About relation of stochastic approach with kinetic equations}

For any particle $x_{i}$ on the real line there is a mark - velocity
$v_{i}>0$. Then the flow of cars is defined by countable subset $\{(x_{i},v_{i})\}$
of the phase space $R\times R_{+}$. And the random flow, at a given
time, is a probability measure $\mu$ on the set of such configurations.
Let $n(A)=n_{\mu}(A)$ be the mean number of particles in bounded
set $A$ of the phase space. If there exists a function $f(x,v)$
such that for any $A$\[
n(A)=\int_{A}f(x,v)dxdv,\]
 then $f$ is called one-particle (correlation) function of the measure
$\mu$.

Although kinetic equations approach to transportation problems is
known long ago \cite{Prigogine} and was studied a lot in physical
papers \cite{Helbing_1,Helbing_2} we think that rigorous mathematical
deduction of kinetic equations (that is equations for $f(t;x,v)$)
from dynamics of individual cars stays open. Here we give only short
comments.

Let at time $t=0$ be given a distribution $\mu$ with one-point function
$f(0,x,v)$. If particles move freely, that is each particle moves
with fixed velocity, different for different cars, then \[
f(t+\delta;x,v)=f(t;x-v\delta,v).\]
 Subtracting $f(t;x,v)$ from both parts of this equality, dividing
by $\delta$ and taking limit $\delta\rightarrow0$, we get\begin{equation}
\frac{\partial f}{\partial t}+v\frac{\partial f}{\partial x}=0.\label{kin}\end{equation}
 This is the simplest kinetic equation without collisions. Its solution
is of course \begin{equation}
f(t;x,v)=f(0;x-vt,v).\label{freeKinetics}\end{equation}

In general kinetic equations should contain collisional term. What
is collision in traffic ? It could be close approaching of cars, that
is the reason for the drivers to change their behavior, for example
to perform overtaking. Also one can consider formation of a new compound
particle - cluster of several cars with short lifetime, which decays
shortly on smaller clusters or separate cars. Moreover, there can
be clusters on various scales - such situation had never been formalized
even in mathematical physics.

Another possible approach is to consider movement of fast cars in
the media of slow cars. Let there exist two types of cars with velocities
$v$ and $v_{1}>v$ correspondingly. Slow cars move independently
and their one-particle function looks like (\ref{freeKinetics}).
In the coordinate system moving with speed $v$, slow cars will become
standing obstacles. Rules of overtaking define velocity $v_{1}(x)$
of fast cars, in the new coordinate system, in the vicinity of separate
obstacle. This can give, under appropriate scaling, some kinetic picture.
But as far as we know there were no efforts to write down rigorously
even such simple models.

There are however papers where Burgers equation is deduced from stochastic
particle picture. However, the corresponding particle process is far
from one-dimensional car flow, see reviews on exclusion processes
\cite{Blythe,Derrida}.

\section{Calculation of mean velocity on the road}

We consider here simplest problems concerning diminishing of the road
capacity because of random fixed obstacles (accidents, repair work)
and because of slow cars. Our goal is to show (solving completely
some model problems) that it is possible to get simple formulas which
allow to understand main reasons why the mean speed can decrease.
Our main assumptions is the homogeneity of the road, exits and gates
and other parameters.

\subsection{Road as one-dimensional queuing network}

The following model is taken from \cite{Kelly}, p. 117. Let we have
infinite road and two types of cars, given by points of the real line,
the points move in the right direction. First type (faster) cars move
with constant speed $v_{1}$, cars of the second type (slow cars)
have constant speed $v_{1}>v_{2}.$

Assume that faster cars initially (at time $t=0$) have Poisson configuration
on $R$ with density $\lambda_{1}$. Slow cars initially are at the
points\[
x_{0}=0<x_{1}<...<x_{n}<...,\]
 where the distances $x_{k}-x_{k-1}$ are identically distributed
with mean $\lambda_{2}^{-1}$ (not necessarily exponentially). Slow
cars move completely independently - other cars do not influence them.
Faster cars <<interact>> with any car when their coordinates coincide.
Namely, they can overtake slow cars. When a faster car catches up
with a slow car, that is their coordinates coincide, then it follows
slow car for some time, that is moves with the speed $v_{2}$. In
exponential time with parameter $\mu$ the fast car overtakes the
slow car and moves again with the speed $v_{1}$. If the slow car
catches up the group of fast cars still following some slow car, then
the overtaking occurs as in a queue, in order how the fast cars follow
the slow car.

Without loss of generality one can put $v_{2}=0$, and the speed of
fast cars put equal to $v=v_{1}-v_{2}$. That is why each slow car
can be considered as a server, where the clients (fast cars) arrive
and are waiting for the service (overtaking), they are served with
intensity $\mu$.

Now the problem can be reduced to the linear queuing system which
we shall describe. There is infinite sequence \[
S_{0}\to...\to S_{k}\to S_{k+1}\to...\]
 of nodes (servers) of two types. Each server $S_{k}$ is a $M/M/1$
type server with FIFO (first-in-first-out) discipline, that is in
the natural order. These servers correspond to slow cars, and clients
correspond to fast cars. For example $S_{0}$ corresponds to the extreme
left slow car. Second letter $M$ means exponential service rate.
First letter $M$ means Poisson arrival stream. Thus the clients to
$S_{0}$ arrive as a stationary Poisson stream with intensity $\lambda_{1}v$.
From elementary queuing theory it has been known, firstly, that if
$\lambda_{1}v<\mu$, then the stationary regime will be established
with the probabilities \[
P_{n}=(1-r)r^{n},r=\frac{\lambda_{1}v}{\mu}.\]
 that the length of the queue equals $n$. Secondly, (Burke theorem),
that in the stationary regime the output flow from $M/M/1$ will be
Poisson with intensity equal to intensity of the input flow, which
is $\lambda_{1}v$ in our case.

After the first node, the output flow, with random but the same for
all clients time shift $\frac{x_{1}-x_{0}}{v}$, arrives to node $S_{1}$,
where also stationary regime will be established.

Let us find mean velocity of fast cars on the interval $(x_{0},x_{N}),N\to\infty$.
We shall assume that the stationary regime has already been installed.
The time along this path is the sum from $N$ overtaking and $N$
paths between slow cars. Mean overtaking time is \[
\sum_{n=0}^{\infty}(1-r)r^{n}\frac{(n+1)}{\mu}=\frac{1}{(1-r)\mu}=\frac{1}{\mu-\lambda_{1}v},\]
 and mean time to catch up the next slow car is\[
\frac{1}{\lambda_{2}v}\]
 Thus, the mean speed of fast cars is \[
v_{mean}=\frac{\lambda_{2}^{-1}}{(\mu-\lambda_{1}v)^{-1}+(\lambda_{2}v)^{-1}}.\]

In the next section we will consider more complicated situation with
more general distributions.

\subsection{Mean speed slow down due to repair works}

On a long road the cars move with constant speed $v$, but encounter
obstacles. The obstacles have small size comparative to distances
between them, and we assume the obstacles are just points. They can
appear on arbitrary part of the road $(x,x+dx)\subset R$ during time
interval $(t,t+dt)\subset R$ with probability $\lambda dxdt$. More
exactly, the pairs (coordinate and time moment of obstacle appearance)
$(x_{j},t_{j})\in R\times R_{+}$ form Poisson random filed $\Pi$
on $R\times R_{+}$with intensity $\lambda$. Another equivalent definition
is that for any time interval $I\subset R$ there is Poisson arrival
stream with intensity $\lambda|I|$, moreover the arriving obstacle
choose the point uniformly on the interval $I$.

Assume that $j-$th obstacle stays on the road for some random time
$\tau_{j}$, after that it is cleared off from the road. Random variables
$\tau_{j}$ are assumed to be i.i.d. with distribution function $Q(t)$
and independent from Poisson random field $\Pi$. Assume finiteness
of $m_{Q}=E\tau_{j}$ and $m_{Q}^{(2)}=E\tau_{j}^{2}$.

Further on we consider two cases. In the first case the bypass is
prohibited and the car stands until the obstacle will be deleted,
after this the car immediately moves with speed $v$. In the second
case, assume the new car arrives to an obstacle and sees $n$ cars,
which have not yet bypassed this obstacle. The new car is allowed
to bypass the obstacle and these $n$ cars simultaneously. This takes
takes some random time which does not depend on $n$. Denote $\eta_{m,i}$
the random bypass time of $m$-th obstacle by $i$-th car. Random
variables $\eta_{m,i}$ are assumed to be i.i.d. with distribution
function $F(u)$. These assumptions are natural for small density
of cars, with small number of cars accumulating after the obstacle.
Below we consider the case of greater load.

First of all we shall find the mean speed of the car. Under above
assumptions the cars do not impede each other, that is why it is sufficient
to consider one car problem. Denote $T(x)$ random time the car passes
distance $x$. We will find the limit of the ratio $\frac{x}{T(x)}$
as $x\to\infty$.

Let $b=\lambda m_{Q}$, and $\zeta$ be a random variable with density
\begin{equation}
h(t)=m_{Q}^{-1}(1-Q(t)).\label{h}\end{equation}
Note that\[
E\zeta=\frac{1}{m_{Q}}\intop_{0}^{\infty}t(1-Q(t))dt=\frac{1}{m_{Q}}\intop_{0}^{\infty}(1-Q(t))d\left(\frac{t^{2}}{2}\right)=\frac{1}{2m_{Q}}\intop_{0}^{\infty}t^{2}dQ(t)=\frac{m_{Q}^{(2)}}{2m_{Q}},\]
 where $m_{Q}^{(2)}$ is the second moment of $Q(t).$

Put $\alpha=\min(\eta,\zeta)$, where equality is in distribution.
Moreover, the random variables $\eta$,$\zeta$ are considered independent
and $\eta$ has distribution function $F(u)$ . Put\[
a=E\alpha.\]

\begin{theorem}

With probability $1$ as $x\to\infty$\begin{equation}
\frac{x}{T(x)}\to\frac{v}{1+abv}.\label{mv}\end{equation}

\end{theorem}

Proof. With no loss of generality assume that the car enters the road
at the point $x=0$ at time $t=0.$ Let $T_{0}(x)$ be the idle time
of the car. Then obviously $T(x)-T_{0}(x)=v^{-1}x$ and\[
\frac{x}{T(x)}=\frac{x}{T(x)-T_{0}(x)+T_{0}(x)}=\frac{1}{v^{-1}+x^{-1}T_{0}(x)}\]
 Thus it is sufficient to find the limit of the ratio $\frac{T_{0}(x)}{x}$
as $x\to\infty$. We want to show that \begin{equation}
T_{0}(x)=\sum_{i=1}^{\pi(x)}\alpha_{i},\label{loss}\end{equation}
 where $\alpha_{i}$ are i.i.d. random variables distributed as $\alpha$,
$\pi(x)$ is a random variable with Poisson distribution with parameter
$bx$. Random variables $\alpha_{i}$ and $\pi(x)$ are assumed to
be independent. The meaning of this formula, that the car on the distance
$x$ will meet $\pi(x)$ obstacles and will loose random time $\alpha_{i}$
$i$-th obstacle.

From (\ref{loss}) and strong law of large numbers it follows easily
that $\frac{T_{0}(x)}{x}\to ab$ a.e. as $x\to\infty$.

Let us prove (\ref{loss}). Introduce marked Poisson point field $\Pi_{1}$
on $R\times R_{+}$ with configuration $(x_{j},t_{j},\tau_{j})$,
that is $\tau_{j}$ is the mark at the point $(x_{j},t_{,j})$. The
following assertion can be found in \cite{Daley}:

\begin{lemma}\label{lm0}

Marked point field $\Pi_{1}$ is equivalent (in probability) to the
Poisson field on $R\times R_{+}^{2}$ with intensity $\lambda dxdtdQ(t)$.

\end{lemma}

On picture 1 the obstacles are presented as horizontal segments. The
coordinates of the initial point define place and time of its appearance
(pair $(x_{j},t_{j})$). The length of the segment is the time of
obstacle stay on the road (mark $\tau_{j}$).

\includegraphics[scale=0.7]{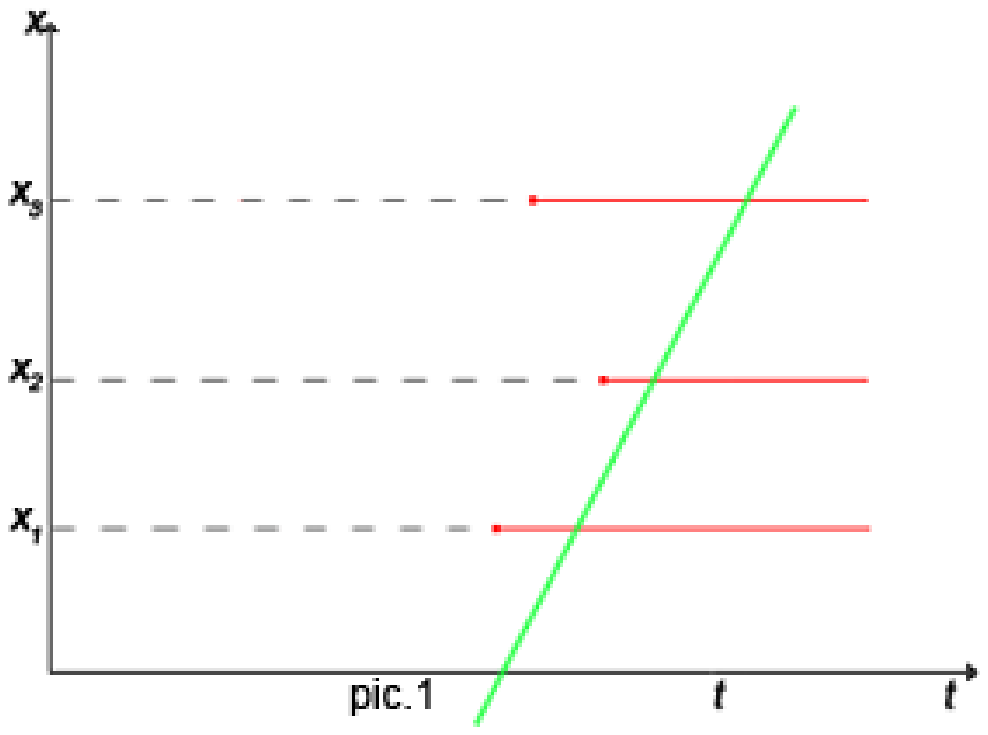}

Consider any straight line $c_{1}t+c_{2}$ and consider its intersection
points with horizontal segments. Denote $\{x_{i}\}$ the space coordinates
of these points as it is shown on picture 1. The next lemma is proved
\cite{CoxIsham}.

\begin{lemma}\label{lm1}

Configuration $\{x_{i}\}$ is the Poisson process with intensity $b=\lambda m_{Q}$.

\end{lemma}

Picture 2 shows the trajectory of the car, which starts from point
$x=0$ at time $t=0$. Denote $x_{i}$ the space coordinates of the
obstacles appearing during movement of this car, $t_{i}$ are the
moments of their appearance, $s_{i}$ -- the moments when the car
meets the obstacles, $u_{i}$ -- time moments when the car gets rid
of them either bypassing or because the obstacle disappear; $\alpha_{i}=u_{i}-s_{i}$
-- car delay $i$-th obstacle.

\includegraphics[scale=0.7]{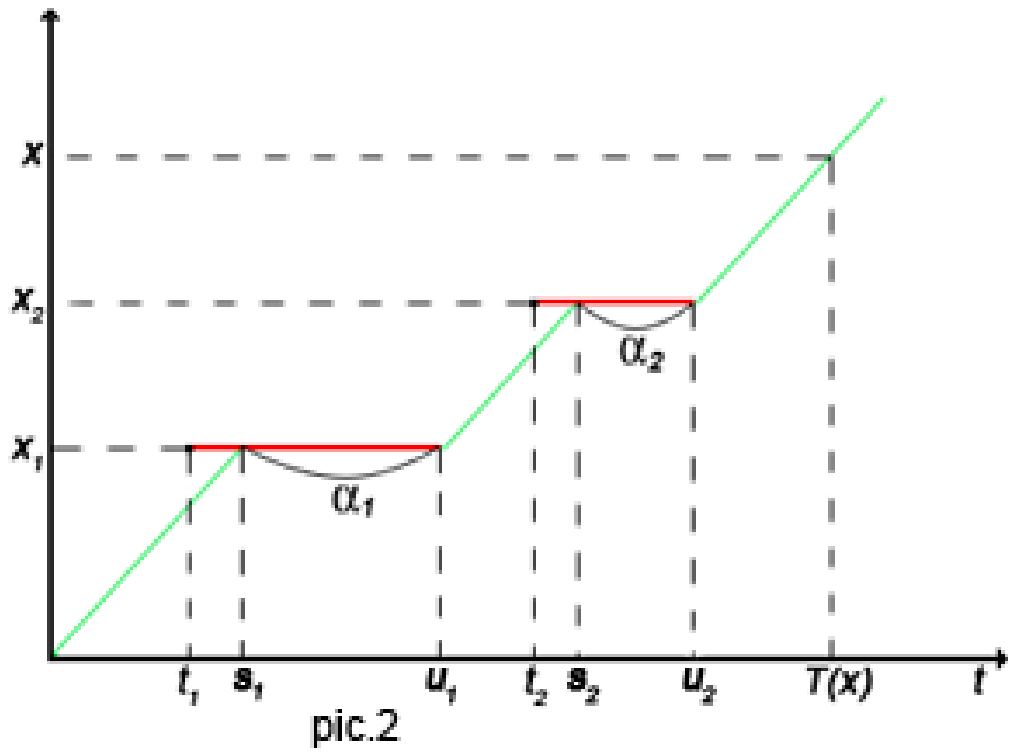}

From lemma \ref{lm1} and from space-time homogeneity of Poisson field
$\Pi$ it follows that the points $x_{i}$ make up Poisson process
of intensity $b$.

The residual lifetime of an obstacle is its lifetime after the car
catches it up. In other words, it is equal to the car delay on this
obstacle.

\begin{lemma}The residual lifetime has the distribution with density
$h(s)$, where $h(s)$ is defined by (\ref{h}).

\end{lemma}

In fact, from the properties of Poisson point field it follows that
the conditional distribution of residual lifetime, under the condition
that the full lifetime equals $t$, coincides with the uniform distribution
on the interval $[0,t]$. By lemma \ref{lm1} the probability of obstacle
appearance on the interval $dx$ equals $\lambda m_{Q}dx+o(dx)$,
and the probability of appearance of an obstacle with full lifetime
$t$ in the interval $dx$ is equal to $\lambda tdQ(t)dx+o(dx)$,
that follows from lemma \ref{lm0}. As \[
\frac{\lambda tdQ(t)dx+o(dx)}{\lambda m_{Q}dx+o(dx)}=\frac{tdQ(t)}{m_{Q}}\]
 is the conditional probability of appearance of an obstacle with
lifetime $t$, then the density of residual lifetime is \[
\intop_{s}^{\infty}\frac{tdQ(t)}{m_{Q}}\frac{ds}{t}=m_{Q}^{-1}(1-Q(s))ds=h(s)ds.\]
 Lemma is proved.

If bypassing is possible, the car will lose time equal to minimum
of bypass time and residual lifetime of the obstacle, that is $\alpha_{i}=\min(\eta,\zeta)$.
Theorem is proved.

Let us discuss the result. We have seen already the meaning of the
constant $a$, the constant $b$ characterizes stationary density
of obstacles in the space.

This result is rather exact for small density of cars, as there will
not more than one car at the obstacle. For large density of cars the
bypass time will increase dependent on the mean length of the queue
at this obstacle.

\subsection{Mean speed slow down due to slow cars}

The road here is the real line $R$. Flows are not too dense, thus
the length of cars does not play role, and at a given moment the position
of the car is given by the point $x_{i}(t)\in R$, where $i$ is the
index enumerating the cars. Each car has predefined route: place and
time of entrance $x_{i,in},t_{i,in}$, and the exit coordinate $x_{i,out}$.
But the exit time $t_{i,out}$ depends on the road capacity. We define
mean velocity of the car $i$ as \[
V_{i}=\frac{x_{i,out}-x_{i,in}}{t_{i,out}-t_{i,in}}.\]
 There are two types of cars - fast and slow, any car moves with constant
speed in the right direction. Fast cars have speed $v_{1}$, slow
cars have speed $v_{2}$, where $v_{1}>v_{2}>0$. Put $v=v_{1}-v_{2}$.
Slow cars move along the road without stops, and fast cars - until
they reach a slow car. After this the fast car $i$ moves with slow
car $j$ for some random time $\tau_{i,j}$, and them bypasses it
immediately getting speed $v_{1}$. Main assumption is that random
variables $\tau_{i,j}$ are i.i.d. with distribution function $F(s)$.

This distribution function may be found statistically by two ways
- by direct sampling and by estimating the density of obstacles, for
example for one-lane road it can be the density of the opposite flow.

Arrivals of slow cars is given by the same Poisson point field $\Pi$
with intensity $\lambda$, which was defined in the previous section.
We will need more notation. To any slow car $j$ we assign random
distance $\rho_{j}$ which it should pass, after that it will leave
the road. Random variables $\rho_{i}$ are assumed i.i.d. with distribution
function $G(r)$. Moreover $\rho_{j}$ does not depend on the random
field $\Pi$. Assume existence of two first moments $m_{G}=E\rho_{1}$,
$m_{G}^{(2)}=E\rho_{1}^{2}$.

Slow car does not meet obstacles and passes all its route with speed
$v_{2}.$ Fast cars can be stopped by slow cars. Consider two cases.
In the first case the overtaking is prohibited and fast car follows
the slow car until the latter will leave the road, after that the
fast car immediately gets its speed $v_{1}$. In the second case the
overtaking is possible. More exactly, when $i$-th fast car catches
up $j$-th slow car or a group of fast cars following $j$-th slow
car, it takes random time $\tau_{i,j}$ to bypass $j$-th slow car
or all group of these cars. This time does not depend on the size
of the group. $\tau_{i,j}$ are assumed to be i.i.d. with distribution
function $F(u)$.

Put $d=\lambda m_{G}\left(v_{2}^{-1}-v_{1}^{-1}\right)$. Introduce
random variable $\beta$ with distribution density $g(x)=m_{G}^{-1}(1-G(x))$
and put $\gamma=\min(v_{2}\tau_{1,1},\beta)$, where the equality
is in distribution, and random variables $\tau_{1,1}$,$\beta$ are
assumed independent. Note that \[
E\beta=\frac{m_{G}^{(2)}}{2m_{G}}.\]
 Put $c=E\gamma$.

\begin{theorem}

With probability $1$ as $x\to\infty$\begin{equation}
\frac{x}{T(x)}\to\bar{v}_{1}=\frac{1+dc}{1+dcv_{1}v_{2}^{-1}}v_{1}.\label{mv2}\end{equation}

\end{theorem}

Proof. Let us note that this case can be reduced to the case $v_{2}=0.$
Introduce new coordinate system moving with the speed $v_{2}$ relative
to the initial system. Find the mean speed of fast cars in the new
coordinate system with the formula (\ref{mv}), substituting $v=v_{1}-v_{2}$,
$b=\frac{\lambda m_{G}}{v_{2}}$, $a=\frac{c}{v_{2}}$: \[
\frac{1}{(v_{1}-v_{2})^{-1}+\frac{\lambda mc}{v_{2}^{2}}}=\frac{v_{1}-v_{2}}{1+dcv_{1}v_{2}^{-1}}.\]
 Then the mean speed of a fast car relative to the initial coordinate
system is\[
\bar{v}_{1}=\frac{v_{1}-v_{2}}{1+dcv_{1}v_{2}^{-1}}+v_{2}=\frac{1+dc}{1+dcv_{1}v_{2}^{-1}}v_{1}.\]

\section{Analysis of complex transport system}

Graph is a natural tool to describe transport networks (for example,
city streets), where the set $V$ of vertices represent crossroads
(nodes or service stations), and the set $L=\{(i,j)\}$ of edges represent
segments of streets without crosses. Let $N$ be the number of vertices.
We assume that there not more than one edge between two vertices.

Mostly the following two classes of networks were studied. First class
includes (by the name of the authors and in the order of increasing
generality) Jackson networks, BCMP and DB networks (see \cite{Kelly},
\cite{Serfozo}). Therein the client (communication, car, job), after
being served in some vertex, chooses randomly the next vertex. The
second class are the Kelly networks (see \cite{Kelly}), where each
client has apriori fixed route. These classes have much in common
- close results and techniques. For example, they may have very useful
multiplicativity property, that is the stationary distribution has
the so called product form. Here we consider only the first class
of networks.

\subsection{Closed networks}

If the cars do not come from without and do not leave outside, the
network is called \textbf{closed}. Then the number of cars in the
network is conserved and we denote it by $M$. The movement of a car
is defined as follows. The car is waiting for some time on the cross
$i$ and is directed afterwards to some vertex $j$. The choice of
$j$ is defined by some stochastic matrix, we shall call it the routing
matrix, $P=\{p_{ij}\}_{i,j=1,...,N}$, where $p_{ij}$ is the probability
that from vertex $i$ the car will go (after service time) to vertex
$j$ (for example, to the right, left or straight on), that is along
the street $(i,j)$.

Stochastic matrix $P$ defines finite discrete time Markov chain with
state space $V=\{1,...,N\}$. We assume this Markov chain irreducible.
Then the system of linear equations\begin{equation}
\rho P=\rho,\:\rho=(\rho_{1},...,\rho_{N})\Longleftrightarrow\sum_{i=1}^{N}\rho_{i}p_{ij}=\rho_{j},\: j=1,...,N\label{traff}\end{equation}
 has a unique solution (up to a common factor). Normed solution gives
its stationary distribution\[
\pi_{i}=\frac{\rho_{i}}{\sum_{i=1}^{N}\rho_{i}},\: i=1,...,N\]
 For each vertex $i\in V$ we define a function $\mu_{i}(n_{i})$
of the number $n_{i}$ of cars in $i$-th vertex, where $\mu_{i}(0)=0$
and $\mu_{i}(n_{i})>0$ for $n_{i}>0$. This function characterizes
capacity of this node and defines the intensity of output flow from
this vertex. Namely, the probability that from node $i$, during time
period $t$, exactly one car will leave the vertex is $\mu_{i}(n_{i})dt+o(dt)$,
under the condition that in the node there are $n_{i}$ cars at time
$t$. Using queuing theory terminology we call $\mu_{i}(n_{i})$ the
service intensity in vertex $i$.

The order in which the arriving clients are served is defined by service
discipline (protocol). The simplest possibility is FIFO - first come
first served. If there are $n_{i}$ cars in the vertex $i$, then
the first car in the queue is served with intensity $\mu_{i}(n_{i})$.

More general service discipline is the resource sharing discipline,
the resource here is the capacity of the crossroad. This means that
resource is divided in some proportion between cars in this crossroad.
We assume that that $k$-th car in $i$- th node is served with intensity
$\mu_{i,k}(n_{i})\le\mu_{i}(n_{i})$ so that \[
\sum_{k=1}^{n_{i}}\mu_{i,k}(n_{i})=\mu_{i}(n_{i}).\]
 For example, common resource can be divided equally between each
car in the queue \[
\mu_{i,k}(n_{i})=\frac{\mu_{i}(n_{i})}{n_{i}}\]
 If $\mu_{i,1}(n_{i})=\mu_{i}(n_{i})$, then we get FIFO discipline.
Thus, the intensities $\mu_{i,k}(n_{i})$ completely define service
discipline in the nodes.

Dynamics in the network is described by $N-$dimensional continuous
time Markov chain $\xi(t)=(\xi_{i}(t),\, i=1,...,N)$, where $\xi_{i}(t)$
is the number of cars-at vertex $i$ at time $t$. Thus, the state
space $S_{M}$ is the set of all vectors with nonnegative integer
coordinates $\bar{n}=(n_{1},...,n_{N})$, such that $n_{1}+...+n_{N}=M$.

Let $e_{i}$ be a base vector, where $i$-th coordinate is 1, and
the rest coordinates are $0$. From state $\bar{n}$ Markov chain
$\xi(t)$ can jump to the state $T_{i,j}\bar{n}=\bar{n}-e_{i}+e_{j},i\neq j$,
with intensity \begin{equation}
\alpha(\bar{n},T_{i,j}\bar{n})=\mu_{i}(n_{i})p_{i,j},\label{trans}\end{equation}
 under the condition that $n_{i}\neq0$. The jump $\bar{n}\rightarrow T_{i,j}\bar{n}$
corresponds to the fact that after leaving cross $i$ the car arrives
to cross $j$.

Note that Markov chain $\xi(t)$ is uniquely defined by the routing
matrix $P$ and service intensities $(\mu_{i}(n_{i}),i=1,...,N)$.

Let $\rho=(\rho_{1},...,\rho_{N})$ be the solution of equations (\ref{traff}),
which we consider as a formal equation for intensities $\rho_{i}$
of input flows to the nodes (in the stationary regime these are equal
to output intensities). Solving these equations, we find $\rho_{i}$.
Then the stationary distribution $\nu(n_{1},...,n_{N})$ of Markov
chain $\xi(t)$ looks like\begin{equation}
\nu(n_{1},...,n_{N})=\frac{1}{Z_{N,M}}\prod_{i=1}^{N}\,\frac{\rho_{i}^{n_{i}}}{\mu_{i}(1)\mu_{i}(2)...\mu_{i}(n_{i})},\label{equil}\end{equation}
 where the normalizing factor (the canonical partition function) \[
Z_{N,M}=\sum_{n_{1}+...+n_{N}=M}\prod_{i=1}^{N}\,\frac{\rho_{i}^{n_{i}}}{\mu_{i}(1)\mu_{i}(2)...\mu_{i}(n_{i})},\]
 This can be verified by direct substitution of (\ref{equil}) to
Kolmogorov equations for stationary probabilities, see for example
\cite{Kelly}.

\subsection{Open networks}

Consider network consisting of $N$ nodes. Contrary to close network,
total number of cars is not fixed. Assume that to node $i$ arrives
(from outside) Poisson flow with intensity $\lambda_{i}$, $i\in\{1,...,N\}$.

The routing matrix $P=\{p_{ij}\}_{i,j=1,...,N}$ is assumed to be
decomposable and \begin{equation}
\forall\, i\,:\:\sum_{j=1}^{N}p_{ij}\leq1,\quad\exists\, i_{0}:\:\sum_{j=1}^{N}p_{i_{0}j}<1.\label{subst}\end{equation}
 Similar to closed network, $p_{i,j}$ is the probability that a car
from node $i$ goes to node $j$. The difference is that there is
also probability\[
p_{i0}=1-\sum_{j=1}^{N}p_{ij}.\]
 for a car at node $i$ to leave the network.

As for closed network, let $\mu_{i}(n_{i})$ be the service intensity
at node $i$. Then the car at node $i$ leaves the network with intensity
$\mu_{i}(n_{i})p_{i,0}$.

The dynamics of an open network is described with $N$-dimensional
continuous time random vectors $\eta(t)=(\eta_{i}(t),\, i=1,...,N)$,
where $\eta_{i}(t)$ is the number of cars at node $i$ at time $t$.
Random process $\eta(t)$ is a continuous time Markov chain with state
space $S$, where $S$ is the set of $N$-dimensional vectors with
non-negative integer components $\bar{n}=(n_{1},...,n_{N})$. From
the state $\bar{n}$ Markov chain $\xi(t)$ can jump to one of the
states $T_{i,j}\bar{n}=\bar{n}-e_{i}+e_{j}$, $T_{i,0}\bar{n}=\bar{n}-e_{i}$,
$T_{i}\bar{n}=\bar{n}+e_{i}$ with intensities\begin{eqnarray}
\alpha(\bar{n},T_{i,j}\bar{n}) & = & \mu_{i}(n_{i})p_{i,j},\nonumber \\
\alpha(\bar{n},T_{i,0}\bar{n}) & = & \mu_{i}(n_{i})p_{i,0},\label{trans1}\\
\alpha(\bar{n},T_{i}\bar{n}) & = & \lambda_{i},\nonumber \end{eqnarray}
 under condition that $T_{i,j}\bar{n},\, T_{i,0}\bar{n},\, T_{i}\bar{n}\in S$.

Thus, Markov chain $\eta(t)$ is uniquely defined by the triplet $(\lambda,\mu,P)$,
where $\lambda=(\lambda_{1},...,\lambda_{N})$ is the vector of intensities
of input flows, $\mu=(\mu_{i}(n_{i}),i=1,...,N)$ -- the vector of
service intensities in the nodes and $P$ is the routing matrix.

Consider formal equation for input intensities to the nodes (in the
stationary regime they are equal to output intensities) \begin{equation}
\rho=\lambda+\rho P\quad\Longleftrightarrow\rho_{i}=\lambda_{i}+\sum_{k=1}^{N}\rho_{k}p_{ki},\;\forall\, i\label{traff_o}\end{equation}
 Under the condition (\ref{subst}) and the irreducibility of $P$,
this equation has a unique solution which can be written as follows\[
\rho=\lambda+\sum_{n=1}^{\infty}\lambda P^{n}.\]

Consider now the case when the service intensities $\mu_{i}(n_{i})\equiv\mu_{i}$
do not depend on the number of cars in nodes. Define the loads in
the nodes by the formula\[
r_{i}=\frac{\rho_{i}}{\mu_{i}},\, i=1,...,N.\]

The following theorem can be found, for example in \cite{Kelly},
\cite{FMM}, it is called Gordon-Newell.

\begin{theorem}Markov chain $\eta(t)$ is ergodic if and only if
for all $i=1,...,N$ we have $r_{i}<1$. Then the stationary distribution
is given by\[
\sigma(n_{1},...,n_{N})=\prod_{i=1}^{N}(1-r_{i})r_{i}^{n_{i}}.\]
 \end{theorem}

It easily follows from this theorem that the mean queue lengths in
the stationary regime are-\[
m_{i}=\frac{r_{i}}{1-r_{i}}.\]
 If in some nodes $i_{1},...,i_{k}$ the load is strictly larger than
1, then Markov chain$\eta(t)$ is transient. Then the mean queue lengths
at nodes $i_{1},...,i_{k}$ tend to infinity as $t\to\infty$. A detailed
analysis of open networks see in \cite{BotZam}. In particular, it
is shown that in the nodes with the load exceeding $1$, mean queue
lengths grow linearly with time. There are also explicit formulas
for the mean speed of growth.

\subsection{Algorithm to find the critical load in closed networks}

This section is based on the paper \cite{MalYak}. Consider a sequence
of closed networks $J_{N}$, $N=1,2,...$. Network $J_{N}$ consists
of $N$ nodes and $M=M(N)$ cars. Service intensities at the nodes
of network $J_{N}$ do not depend on the queue length : $\mu_{i,N}(n_{i})\equiv\mu_{i,N}$.
Let $P_{N}=\{p_{i,j,N}\}$ be the routing matrix in the $N$-th network;
$P_{N}$ is assumed to be irreducible.

Let $\rho_{N}=(\rho_{1,N},...,\rho_{N,N})$ be the vector with positive
components satisfying equation \begin{equation}
\rho_{N}=\rho_{N}P_{N}.\label{traff1}\end{equation}
 Relative loads in the nodes are defined as\[
r_{i,N}=C_{N}^{-1}\rho_{i,N}\tau_{i,N},\]
 where $\tau_{i,N}=\mu_{i,N}^{-1}$ and $C_{N}=\max_{i=1,...,N}\rho_{i,N}\tau_{i,N}$.
It is evident that $r_{i,N}\in[0,1]$.

Correspondingly to (\ref{equil}) the stationary distribution of the
number of cars $\xi_{i,N,M}$ in the nodes of $J_{N}$ is equal to\[
P_{N,M}(\xi_{i,N,M}=n_{i},i=1,...,N)=\frac{1}{Z_{N,M}}\prod_{i=1}^{N}r_{i,N}^{n_{i}},\]
 where the normalizing factor (canonical partition function) \begin{equation}
Z_{N,M}=\sum_{n_{1}+...+n_{N}=M}\prod_{i=1}^{N}r_{i,N}^{n_{i}}.\label{sts}\end{equation}
 Many important network characteristics can be expressed in terms
of this partition function.

\begin{exer}

Show that the mean number of cars in the $i$-the node in the stationary
regime is\begin{equation}
m_{i,N,M}=E\xi_{i,N,M}=\frac{r_{i,N}}{Z_{N,M}}\frac{\partial Z_{N,M}}{\partial r_{i,N}}\label{min-1}\end{equation}

\end{exer}

We are interested in cases when $N,M$ are large enough, more exactly
$N,M\to\infty$ and so that $\frac{M}{N}\to\lambda=const$, that is
the number of cars per one node is constant. We shall see that $\lambda$
defines existence or non-existence of jams in the network.

To formulate results as general as possible, we shall assume weak
convergence of relative loads $r_{i,N}$. More exactly, define sample
measure on the interval $[0,1]$\[
I_{N}(A)=\frac{1}{N}\sum_{i:r_{i,N}\in A}1,\]
 where $A$ is arbitrary Borel subset of $[0,1]$. Assume that as
$N\to\infty$ the measures $I_{N}$ weakly converge to some probability
measure $I$ on $[0,1]$.

\begin{remark}

It could be interesting to find concrete sequences of growing networks
for which the limiting measure $I$ can be found explicitly. Some
examples where measure $I$ is one-point one can find in the bibliography
to the paper \cite{MalYak}.

\end{remark}

In terms of limiting measure $I$ we shall find critical value $\lambda_{cr}$
of the density, so that for $\lambda<\lambda_{cr}$ mean lengths of
the queues are uniformly bounded. If $\lambda\geq\lambda_{cr}$, then
at the node with maximal value of load the mean length of the queue
tends to infinity, that means the existence of a jam at this node.
Put\[
h(z)=\int_{0}^{1}\frac{r}{1-zr}dI(r),\]
 where $z\in\mathbb{C}\setminus[1,+\infty)$. The function $h(z)$
is strictly increasing on $[0,1)$. Put\[
\lambda_{cr}=\lim_{z\to1-}h(z).\]
 We shall assume that $\lambda_{cr}>0$.

\begin{theorem}\label{th1} 
\begin{itemize}
\item If $\lambda<\lambda_{cr}$, then mean lengths of queues are uniformly
bounded: there exists constant $B$, such that $m_{i,N}<B$, uniformly
in $N\geq1$ and $1\leq i\leq N$. 
\item If $\lambda\geq\lambda_{cr}$ and $i(N)$ satisfies condition $r_{i(N),N}=1$,
then $m_{i(N),N}\to\infty,$ as $N\to\infty$, that is jams will be
in the nodes where the load is maximal. 
\end{itemize}
\end{theorem}

For $z\in\mathbb{C}\setminus[1,+\infty)$ put \begin{equation}
S_{N}(z)=-\lambda(1+\varepsilon_{N})\ln z-\frac{1}{N}\sum_{i=1}^{N}\ln(1-zr_{i,N}),\label{snz}\end{equation}
 \[
S(z)=-\lambda\ln z-\int_{0}^{1}\ln(1-zr)dI(r),\]
 where $\lambda(1+\varepsilon_{N})=\frac{M}{N}$.

Introduce the generating function (the grand partition function)\[
\Xi_{N}(z)=\sum_{M=0}^{\infty}z^{M}Z_{N,M}=\prod_{i=1}^{N}\frac{1}{1-zr_{i}},\:|z|<1.\]
 By Cauchy formula and (\ref{snz}) we have the following expression
for the partition function (\ref{sts}):\begin{equation}
Z_{N,M}=\frac{1}{2\pi i}\intop_{\gamma}\,\frac{\Xi_{N}(z)}{z^{M+1}}\, dz=\frac{1}{2\pi i}\intop_{\gamma}\,\frac{\exp(NS_{N}(z))}{z}\, dz,\label{sts1}\end{equation}
 where $\gamma=\{z\in\mathbb{C}:|z|=\sigma<1\}$. For the means, accordingly
to (\ref{min-1}), we have\begin{equation}
m_{i,N}=\frac{1}{2\pi iZ_{N}}\intop_{\gamma}\,\frac{r_{i,N}}{1-zr_{i,N}}\exp(NS_{N}(z))\, dz.\label{min1}\end{equation}
 One can show that for the stationary distribution of queue lengths
the following formula holds\begin{equation}
P_{N,M}(\xi_{1,N,M}=n_{1},...,\xi_{K,N,M}=n_{K})=\frac{1}{2\pi iZ_{N}}\intop_{\gamma}\, z^{-1}\prod_{i=1}^{K}(1-zr_{i,N})(zr_{i,N})^{n_{i}}\exp(NS_{N}(z))dz.\label{pn}\end{equation}

In the proof of theorems of this section, essential role is played
by the steepest descent method (see, \cite{Fedor}), more exactly,
its generalization concerning the case when the function in the exponent
depend on $N$. From equation \begin{equation}
\frac{\partial S_{N}(z)}{\partial z}=0\label{sn}\end{equation}
 we find steepest descent points. Let $z_{0,N}$ be the root of this
equation lying in $(0,1).$

\begin{exer}

Show that all roots of equation (\ref{sn}) are real and positive.
Always there exists unique root, which belongs to the interval $(0,1).$

\end{exer}

Let $z_{0}$-be the root of equation \begin{equation}
h(z)=\frac{\lambda}{z}\:\Leftrightarrow\:\frac{\partial S(z)}{\partial z}=0,\label{eq-h}\end{equation}
 lying in the interval $(0,1).$

\begin{exer} 
\begin{itemize}
\item Prove that for all $\lambda$ there exists limit $\lim_{N\to\infty}z_{0,N}=z_{0}=z_{0}(\lambda)>0$. 
\item If $\lambda<\lambda_{cr}$, then $z_{0}(\lambda)$ is the root of
equation (\ref{eq-h}); $z_{0}(\lambda)$ is strictly increasing in
$\lambda,$ $z_{0}(\lambda)\in(0,1),$ $\lim_{\lambda\to\lambda_{cr-}}z_{0}(\lambda)=1.$ 
\item If $\lambda\ge\lambda_{cr}$, then $z_{0}=1$. 
\end{itemize}
\end{exer}

In the next theorem we find the asymptotics of the partition function
and limiting distribution for the sequence of closed networks $J_{N}$.

\begin{theorem}\label{th2}

Let $\lambda<\lambda_{cr}$, then 
\begin{itemize}
\item As $N\to\infty$ the partition function $Z_{N}$ and the free energy
$F_{N}=\frac{1}{N}\ln Z_{N}$ have the following asymptotics\[
Z_{N}\sim\frac{\exp(NS_{N}(z_{0,N}))}{z_{0}\sqrt{2\pi NS^{\prime\prime}(z_{0})}},\; F_{N}=\frac{1}{N}\ln Z_{N}\sim S(z_{0}).\]

\item If for $i=1,...,K$ there exist the limits $r_{i}=\lim_{N\to\infty}r_{i,N}$,
then\[
\lim_{N\to\infty}m_{i,N}=\frac{z_{0}r_{i}}{1-z_{0}r_{i}},\]
 \[
\lim_{N\to\infty}P_{N,M}(\xi_{1,N,M}=n_{1},...,\xi_{K,N,M}=n_{K})=\prod_{i=1}^{K}(1-z_{0}r_{i})(z_{0}r_{i})^{n_{i}}.\]

\end{itemize}
\end{theorem}

Thus in the limit we get an open network consisting of independent
queues.

\paragraph{Proof of theorems \ref{th1} and \ref{th2}}

We give more general result from which the theorems \ref{th1} and
\ref{th2} follow. Let $U_{d}(v)=\{z\in\mathbb{C}:|z-v|<d\}$. Consider
the contour $\gamma=\{z\in\mathbb{C}:|z|=z_{0}(\lambda)\}$.

\begin{theorem}\label{main}

Let $\lambda<\lambda_{cr}$ and $f(\theta,z)$, $\theta\in\Theta,$
be a family of functions holomorphic in the ring $\{z\in C:z_{0}(\lambda)-\delta_{0}<|z|<z_{0}(\lambda)+\delta_{0}\}$
for some $\delta_{0}>0$, uniformly bounded in this ring and such
that for sufficiently small $\epsilon>0$ there exists such $\delta_{u}>0$
and nonzero constant $f_{u}$, such that $|f(\theta,z)/f_{u}-1|<\epsilon$
for $z\in U_{2\delta_{u}}(z_{0})$, $\theta\in\Theta$.

Then, for $N$ sufficiently large, uniformly in $\theta\in\Theta$\[
\frac{1}{2\pi i}\intop_{\gamma}f(\theta,z)\exp(NS_{N}(z))dz=\frac{f_{u}\exp(NS_{N}(z_{0,N}))}{\sqrt{2\pi NS^{\prime\prime}(z_{0})}}(1+\zeta_{N}),\]
 where $|\zeta_{N}|<25\epsilon.$

\end{theorem}

Proof of this theorem is based on the saddle-point method, see \cite{Fedor}.
Difference from the standard situation is that the function in the
exponent depends on $N$. Detailed proof one can find in the original
paper \cite{MalYak}.

Proof of theorem \ref{th2}. Using theorem \ref{main} we prove the
first part of the theorem \ref{th2}. Using (\ref{sts1}) we have\[
Z_{N,M}=\frac{1}{2\pi i}\intop_{\gamma}\,\frac{\exp(NS_{N}(z))}{z}\, dz,\]
 where $\gamma=\{z\in\mathbb{C}:|z|=z_{0}(\lambda)\}$. Putting $f(\theta,z)=z^{-1}$,
$f_{u}=z_{0}^{-1}$ and using theorem \ref{main}, we get that for
any sufficiently small $\epsilon>0$ for sufficiently large $N$

\begin{equation}
Z_{N}=\frac{\exp(NS_{N}(z_{0,N}))}{z_{0}\sqrt{2\pi NS^{''}(z_{0})}}(1+\zeta_{N}),\:|\zeta_{N}|<25\epsilon.\label{partfu}\end{equation}

The second part of theorem \ref{th2} can be proved similarly by using
(\ref{min1}) for the mean queue length and (\ref{pn}) for joint
distribution of the queue lengths.

\begin{exer}Prove third part of the theorem \ref{th2}, using theorem
\ref{main} and formulas (\ref{min1}), (\ref{pn}).

\end{exer}

Proof of theorem \ref{th1}. To prove the first part of theorem \ref{th1},
consider the family of functions\[
f(\theta,z)=\frac{A}{z}+\frac{\theta}{1-z\theta},\:\theta\in\Theta=[0,1],\ A>0,\ f_{u}=\frac{A}{z_{0}}.\]
 Fix small $\epsilon>0$ and choose $\sigma_{u}=\frac{\epsilon}{8}$,
$A=\frac{16z_{0}}{(1-z_{0})\epsilon}$. By theorem \ref{main} we
have for sufficiently large $N$ and all $\theta\in\Theta$\[
\frac{1}{2\pi i}\int\limits _{\gamma}\left(\frac{A}{z}+\frac{\theta}{1-z\theta}\right)\exp(NS_{N}(z))\, dz=\frac{A\exp(NS_{N}(z_{0,N}))}{z_{0}\sqrt{2\pi NS^{''}(z_{0})}}(1+\zeta_{N}),\:|\zeta_{N}|<25\epsilon.\]
 Dividing by $Z_{N}$ and applying (\ref{partfu}) to the right-hand
part of the resulting equality, we get for sufficiently large $N$\[
A+\frac{1}{Z_{N}}\frac{1}{2\pi i}\int\limits _{\gamma}\frac{\theta}{1-z\theta}\exp(NS_{N}(z))\, dz=A\,(1+\zeta_{N}^{'}),\:|\zeta_{N}^{'}|<30\epsilon.\]
 From last equality and formula (\ref{min1}) we have uniform boundedness
of $m_{i,N}$.

Prove the second part of theorem \ref{th1}. For this we need the
following monotonicity property: for any $M_{2}\geq M_{1}>0$ and
any $N\geq1$ the inequality $m_{i,M_{2},N}\geq m_{i,M_{1},N}$ holds.

As $z_{0}(\lambda)$ is strictly increasing in $\lambda,$ $z_{0}(\lambda)\in(0,1)$
and $\lim_{\lambda\to\lambda_{cr}-}z_{0}(\lambda)=1$, then the function
\[
\frac{z_{0}(\lambda)}{1-z_{0}(\lambda)}\]
 is monotone increasing and tends to $\infty$, when $\lambda\nearrow\lambda_{cr}$.
That is why for any $m>0$ there exists such $\lambda'=\lambda'(m)<\lambda_{cr}$,
that \[
\frac{z_{0}(\lambda')}{1-z_{0}(\lambda')}=m+1.\]

Without loss of generality we can assume that $i(N)\equiv1$ and $r_{1,N}=1$.
If we put $M'(N)=[\lambda'N]$, then by theorem \ref{th2} \[
\lim_{N\to\infty}\, m_{1,M'(N),N}=\frac{z_{0}(\lambda')}{1-z_{0}(\lambda')}.\]
 It follows that for sufficiently large $N$\[
m_{1,M'(N),N}>\frac{z_{0}(\lambda')}{1-z_{0}(\lambda')}-1=m.\]
 But $M/N\to\lambda\ge\lambda_{cr}>\lambda'$, that is why for sufficiently
large $N$ we have $M(N)\ge M'(N)$. By monotonicity $m_{1,N}=m_{1,M(N),N}\ge m_{1,M',N}>m$
for sufficiently large $N$. This proves that $m_{1,N}\to\infty$.

\paragraph{Technical generalizations and mathematical problems}

We assumed above instantaneous displacement between crosses. I.e.
we did not take into account time spent along the streets. This can
be easily overpassed by introducing more extended graph. Namely, introduce
additional vertices $u_{ij}$, corresponding to the street between
crosses $i$ and $j$, and mean duration $\tau_{ij}\text{=}\mu_{ij}^{-1}$
of movement along the streets. In the queuing terms this means that
the streets are considered as new service nodes with infinite number
of servers and exponential service time with mean $\tau_{ij}\text{=}\mu_{ij}^{-1}$.

Note that the results of section 3.1 can be generalized to the case,
when the network contains nodes with infinite number of servers. Let,
for example, the network contain one such node ($i=0$) and $\mu_{0,N}(n)=n\nu_{N}$
is the service intensity in this node. Let $\rho_{N}=(\rho_{0,N},...,\rho_{N,N})$
be the solution of equation (\ref{traff1}). Then relative loads can
be defined by formula \[
r_{i,N}=\frac{\mu_{0,N}}{\rho_{0,N}}\frac{\rho_{i,N}}{\mu_{i,N}},\]
 so that $r_{0,N}=1$. According to (\ref{equil}), the stationary
distribution of the queue lengths is\[
P_{N,M}(\xi_{i,N,M}=n_{i},i=1,...,N)=\frac{1}{\widehat{Z}{}_{N,M}}\frac{1}{(M-\sum_{i=1}^{M}n_{i})!}\prod_{i=1}^{N}r_{i,N}^{n_{i}},\]
 where\[
\widehat{Z}{}_{N,M}=\sum_{n_{1}+...+n_{N}\le M}\frac{1}{(M-\sum_{i=1}^{M}n_{i})!}\prod_{i=1}^{N}r_{i,N}^{n_{i}},\]
 and the grand partition function is\[
\widehat{\Xi}_{N}(z)=e^{z}\prod_{i=1}^{N}\frac{1}{1-zr_{i,N}}.\]
 Put\[
q_{i,N}=\frac{r_{i,N}}{p_{N}},\; w=zp_{N},\quad p_{N}=\max_{1\le i\le N}r_{i,N}.\]
 Then\[
\widehat{\Xi}_{N}(w)=e^{w/p_{N}}\prod_{i=1}^{N}\frac{1}{1-wq_{i,N}}.\]
 Under assumption that $p_{N}N\to\alpha>0$ as $N\to\infty$ we can
find critical value of density $\lambda$ with the formula \[
\lambda_{cr}=\alpha^{-1}+\lim_{w\to1-}\intop_{0}^{1}\frac{q}{1-wq}dI(q),\]
 where, as earlier, measure $I$ is the weak limit of sampling measures
for $N\to\infty$ \[
I_{N}(A)=\frac{1}{N}\sum_{i:q_{i,N}\in A}1,\]
 where $A$ is an arbitrary Borel subset of $[0,1]$.

In \cite{FL} for closed networks similar results are obtained for
more general dependence of intensities of the queue lengths.

The trick with the graph extension allows get rid of other restriction:
that for given cross the mean duration of red light is the same for
all directions. It is necessary then, instead of vertex $i$, corresponding
to this cross, introduce several vertices $(i,d)$, where $d$ enumerates
possible directions on cross $i$. This poses some restrictions on
the service times $\tau_{i,d}$ in the new vertices, like\[
\sum_{d}\tau_{i,d}=\tau_{i}.\]

We limited ourselves to the problem, when in the system at least one
jam appear. It is interesting to find number of jams and the mean
number of cars standing in jams.

\paragraph{Relation with real life}

This model is convenient because all parameters can be statistically
estimated. Namely, statistical estimates of the parameters $p_{ij},\mu_{i}$
look like (for example for constant $\mu_{i}$)\[
p_{ij}\sim\frac{N_{ij}(T)}{\sum_{j}N_{ij}(T)},\mu_{i}=\frac{1}{T}\sum_{j}N_{ij}(T),\]
 where $N_{ij}(T)$ is the number of cars in the time interval $[0,T]$,
choosing direction $j$ on the cross $i$.

Practically interesting is the optimization of traffic lights, that
can be achieved by choice of $\tau_{i,d}$ and by traffic lights synchronization
(that is totally ignored in any exponential model), and by changing
matrix $P$ by helping the choice of route.

\end{document}